\documentclass[12pt]{article}
\usepackage{WideMargins2}

\def\Z{{{\mathbb Z}}}

\def\F{{{\mathbb F}}}
\def\KK{{{\mathbb K}}}

\def\cj#1{{{\overline{#1}}}}

\def\deg{{{\rm deg}}}
\def\ang#1{{{\langle #1 \rangle}}}

\let\x=\times

\let\z=\zeta
\def\ie{{ \it i.e.}}

\title{\Large\bfseries  An explanation for the similar appearance
of two identities involving Dickson polynomials}
\author{Antonia W.~Bluher\\ National Security Agency \\ tbluher@access4less.net}
\date{February 20, 2022}
\begin {document}
\fancytitle

\begin{abstract}
This articles explains the similar appearance of two polynomial
identities involving Dickson polynomials in char.~2, one found by Abhyankar, 
Cohen and Zieve, and the other found by the author.
\end{abstract}

\section{Introduction} \label{sec:Intro}

The $k$th Dickson polynomial $D_k(x)\in \Z[x]$ is defined by
the recursion $D_k(x)=xD_{k-1}(x)-D_{k-2}(x)$ for $k\ge 2$, 
with the initial conditions $D_0(x)=2$, $D_1(x)=x$.
It satisfies the functional equation
$D_k(u+1/u)=u^k+u^{-k}$, which could serve as an alternate definition.
If $k\ge1$, then $D_k$ is monic of degree~$k$.

In \cite{ACZ},
Abhyankar, Cohen, and Zieve found an identity of Dickson polynomials
in finite characteristic, which in char.~2 may be written as
follows, with $q=2^n$:
\begin{equation} X^{q^2-1} + (D_{q+1}(Y)/Y) X^{q-1} + Y^{q-1}
=(X^{2q-2} + Y X^{q-1} + 1)\left(
\prod_{w\in\F_q^\x} (D_{q-1}(wX) - Y) \right). \label{aczIdentity}
\end{equation}
In \cite{DicksonIdentity}, the author found an identity
that is special to char.~2:
\begin{equation}
	X^{q^2-1} + (D_{q-1}(Y)/Y) X^{q-1} + Y^{q-1}
= \prod_{w\in\F_q^\x} \left(D_{q+1}(wX) - Y \right).
\label{bluherIdentity}
\end{equation}
It seems as though $D_{q+1}$ and $D_{q-1}$ switch roles in these equations.
The identities were found independently, in different contexts, and with
different applications, yet their visual similarity
begs for an explanation. This article provides such an explanation.

We say a word about how the identities were found.
The authors of \cite{ACZ} were motivated by some results in group
theory which, when combined with Galois theory, led them to seek
bivariate polynomials of a particular form with a particular type of
factorization.  This search led to them to discover the
identity~(\ref{aczIdentity}). 
The discovery of the second identity arose from the author's attempt to 
understand why certain pairs of polynomials in $\F_{2^n}[x]$
have related factorizations.
The identity explained these related factorizations and led to the result that
if $\F$ is any field of char.~2, $0\ne a\in\F$, and $q=2^n>2$, then
$x^{q+1}+x+1/a$ and $C(x)+a$ have the same splitting field over~$\F$, where
$C(x) = x (\sum_{i=0}^{n-1} x^{2^i-1})^{q+1}$ is a
M\"uller--Cohen--Matthews polynomial of degree $(q^2-q)/2$.

Define $\ang u = u + 1/u$;  then the Dickson relation may be written as
$D_k(\ang u) = \ang{u^k}$. Note that 
$$D_k\circ D_\ell(\ang u) =
D_k(\ang{u^\ell}) = \ang{u^{\ell k}} = D_{k\ell}(\ang u),$$
\begin{equation} \ang u \ang v = \ang{uv} + \ang{u/v}.
	\label{wellKnown1}
\end{equation}
These imply the well-known formulas:
\begin{equation*}
	D_k \circ D_\ell(x) = D_\ell \circ D_k(x) = D_{k\ell}(x),\qquad
D_k(x) D_\ell(x) = D_{k+\ell}(x) + D_{|k-\ell|}(x).
\end{equation*}

Let  $X$ be a transcendental over a field of char.~2, $q=2^n$, and
$$  v = D_{q-1}(X),\qquad y = D_{q+1}(X),\qquad z = D_{q^2-1}(X).$$
Let $U$ be a solution to $U^2+UX+1=0$, so $U$ is transcendental and 
\begin{equation} X=\ang U,\qquad  
v = \ang{U^{q-1}}, \qquad 
y= \ang{U^{q+1}},\qquad 
	z = \ang{U^{q^2-1}}. \label{eq:U}
\end{equation}
By (\ref{wellKnown1}), $\ang{U^q}\ang{U}=\ang{U^{q+1}}+\ang{U^{q-1}}$.
Since $\ang{U^q}=\ang{U}^q$ in char.~2, 
\begin{equation}  y+v=X^{q+1}. \label{eq:yv} \end{equation}

The right side of (\ref{aczIdentity}) vanishes when one specializes
$Y=v$, and so 
(\ref{aczIdentity}) implies
\begin{equation*} X^{q^2-1} + (D_{q+1}(v)/v) X^{q-1} + v^{q-1} = 0.
\end{equation*}
Likewise, the right side of (\ref{bluherIdentity}) vanishes
when one specializes $Y=y$, and so (\ref{bluherIdentity}) implies 
$$ X^{q^2-1} + (D_{q-1}(y)/y) X^{q-1} + y^{q-1} = 0.$$
Now $D_{q+1}(v)=D_{q+1}\circ D_{q-1}(X)=D_{q^2-1}(X)=z$,
and similarly $D_{q-1}(y) = z$.  Thus, the above two formulas may be
written as
$$ X^{q^2-1} + (z/v) X^{q-1} + v^{q-1}=0,\qquad
 X^{q^2-1} + (z/y) X^{q-1} + y^{q-1}=0.$$
Equivalently,
\begin{equation}
	vX^{q^2} + z X^q + v^q X = 0, \label{eq:v}
\end{equation}
\begin{equation}
	yX^{q^2} + z X^q + y^q X = 0. \label{eq:y}
\end{equation}

We showed that (\ref{aczIdentity}) implies (\ref{eq:v}) and
(\ref{bluherIdentity}) implies (\ref{eq:y}), but in fact (\ref{eq:v})
and (\ref{eq:y}) are easy to prove directly. 
Indeed, (\ref{eq:v}) follows from the calculation:
\begin{eqnarray*}
	vX^{q^2}  + z X^q + v^q X &=&  \ang{U^{q-1}} \ang{U^{q^2}} +
\ang{U^{q^2-1}} \ang{U^q} + \ang{U^{q(q-1)}}\ang U \\
&=& (\ang{U^{q^2+q-1}} + \ang{U^{q^2-q+1}})  + 
	(\ang{U^{q^2+q-1}} + \ang{U^{q^2-q-1}}) +  \\
	&& \qquad (\ang{U^{q^2-q+1}} + \ang{U^{q^2-q-1}}) \qquad\qquad\qquad
	\quad\qquad\qquad \text{by~(\ref{wellKnown1})}\\
&=& 0 
\end{eqnarray*}
and (\ref{eq:y}) can be proved similarly. 
The sum of
(\ref{eq:v}) and (\ref{eq:y}) is $(v+y)X^{q^2}+(v+y)^qX=0$,
which is consistent with the well-known formula (\ref{eq:yv}).

\section{Derivation of the identities}
We showed that (\ref{aczIdentity}) easily implies (\ref{eq:v}),
and (\ref{eq:v}) has a one-line proof. Similarly, 
(\ref{bluherIdentity}) easily implies (\ref{eq:y}), and
(\ref{eq:y}) has a one-line proof. Finally, we showed that
(\ref{eq:v}) and (\ref{eq:y}) are related by the well-known
formula (\ref{eq:yv}).  To complete the circle of ideas, we will
show that (\ref{aczIdentity}) can be derived from the known formula (\ref{eq:v}) and 
(\ref{bluherIdentity}) can be derived from the known formula (\ref{eq:y}).
Much of the reasoning given in this section, as well as the one-line
proofs of (\ref{eq:v}) and (\ref{eq:y}) from the previous section, 
can be found in the original proofs of the two polynomial identities; 
see \cite{ACZ,DicksonIdentity}.  The new idea is simply
to show how the similar appearance of (\ref{eq:v}) and (\ref{eq:y}) gives rise
to the similar appearance of the two identities.

For most of the article, $q$ denotes a power of~2, but in the next lemma
$q$ can be any prime power. If $\F$ is a field, then $\cj\F$ denotes its
algebraic closure.

\begin{lemma} \label{lem:fwx}
	Let $\KK$ be a field containing~$\F_q$,
	$f,g \in \KK[X]$, and $F(X)=f(X^{q-1})$. Suppose that
$G(X)=\prod_{w\in\F_q^\x} g(wX)$ has no repeated roots in~$\cj\KK$.
If every root of $g$ in $\cj\KK$ is also a root of $F$, then 
	$G$ divides $F$.
\end{lemma}
\begin{proof}  
	Let $r\in\cj\KK$ be any root of $G$. Then 
	$g(wr)=0$ for some $w\in\F_q^\x$, \ie, $s=wr$ is a root of $g$. 
	By hypothesis, $F(s)=0$. Then $F(r)=f((s/w)^{q-1})=f(s^{q-1})=F(s)=0$.
	We have shown that every root
	of $G$ in $\cj\KK$ is also a root of $F(X)$.  Since $G$ has
	no repeated roots, $G$ divides $F$.
\end{proof}

\begin{lemma} \label{lem:distinctRoots} 
	Let $\KK= \F_q(Y)$, where $Y$ is transcendental and $q=2^n$. If $k>0$
	is odd, then the polynomial 
	$G(X)=\prod_{w\in\F_q^\x} (D_k(wX)-Y)\in \KK[X]$ 
	has no repeated roots in $\cj\KK$.
\end{lemma}
\begin{proof} Let $x\in \cj \KK$ be a root of $D_k(X)-Y$. 
	Write $x=\ang{u}$, where $u\in\cj \KK$. Then $Y=\ang{u^k}$.
	Since $Y$ is transcendental over $\cj\F_q$, so is $u$.
	Let $\mu_k$ denote the $k$th roots of unity in $\cj \F_q$.
	Note that $|\mu_k|=k$ since $k$ is prime to the characteristic.
	If $\z\in\mu_k$, then $D_k(\ang{\z u})=\ang{(\z u)^k} = \ang{u^k} = Y$.
	Thus, $\ang{\zeta u}$ for $\z\in\mu_k$ are roots of $D_k(X)-Y$,
	and 
	$$\text{\it $\{\,w\ang{\zeta u}: w\in\F_q^\x, \z \in \mu_k\,\}$ 
	are roots of $G(X)$.}$$

	We claim these are distinct. To see this,
	suppose $w_1\ang{\zeta_1u} = w_2 \ang{\zeta_2u}$. Then
	$(w_1\zeta_1+w_2\zeta_2)u + (w_1/\zeta_1+w_2/\zeta_2)u^{-1}=0$.
	Since $u$ is transcendental, both coefficients are zero, therefore
	$w_1\zeta_1+w_2\zeta_2 = 0$, $w_1 \zeta_2 + w_2 \zeta_1=0$.
	On summing, we find $(w_1+w_2)(\zeta_1+\zeta_2)=0$,
	so $w_1=w_2$ or $\zeta_1=\zeta_2$.
	If $w_1=w_2$, then the equation $w_1\zeta_1+w_2\zeta_2=0$ implies
	$\zeta_1=\zeta_2$.  If $\zeta_1=\zeta_2$, then the same equation shows
	$w_1=w_2$.  This proves the claim.
	Since there are $(q-1)k$ distinct roots 
	$w\ang{\zeta u}$ and $\deg_X(G)=(q-1)k$, they account for all the
	roots of~$G$, therefore $G$ has distinct roots.
\end{proof}

\noindent{\it Proof of the identity~(\ref{bluherIdentity}).}\quad  Set
$F(X) = X^{q^2-1} + (D_{q-1}(Y)/Y)X^{q-1}+ Y^{q-1}$ and 
$g(X)=D_{q+1}(X)-Y$, considered as polynomials in $\KK[X]$ where
$\KK=\F_{q}(Y)$. Let $x\in\cj\KK$ be any root of $g$ and write $x=\ang u$,
where $u\in\cj\KK$.
Then $Y=\ang{u^{q+1}}$. Since $Y$ is transcendental, so is~$u$.
In (\ref{eq:U}), replace $U$ by $u$; then
$X=\ang{U}$, $y=\ang{U^{q+1}}$, and $z=\ang{U^{q^2-1}}$ are replaced by
$x=\ang{u}$, $Y=\ang{u^{q+1}}$, and $D_{q-1}(Y)=\ang{u^{q^2-1}}$.
Eq.~(\ref{eq:y}) becomes $YxF(x)=0$, so $F(x) = 0$. This shows that 
every root of $g$ is a root of $F(X)$.  
By Lemmas~\ref{lem:fwx} and~\ref{lem:distinctRoots},
$\prod_{w\in\F_q^\x}g(wX)$ divides $F(X)$.
Both polynomials are monic in $X$ of degree~$q^2-1$ and one divides the other,
so they are equal. This proves~(\ref{bluherIdentity}). \qed
\medskip

\noindent{\it Proof of the identity (\ref{aczIdentity}).}\quad
Let $X$ and $V$ be independent transcendentals over $\F_q$.
We will prove the identity (\ref{aczIdentity}) with $Y$ replaced by $V$. 
Let $\KK=\F_q(V)$ and define
$$F(X) = X^{q^2-1} + (D_{q+1}(V)/V)X^{q-1} + V^{q-1},\qquad
G(X)=\prod_{w\in\F_q^\x}(D_{q-1}(wX)-V)$$
considered as elements of $\KK[X]$.  Let $x\in\cj\KK$ be a root of 
$D_{q-1}(X)-V$ and write $x=\ang u$, where $u\in\cj\KK$.
Then $V=\ang{u^{q-1}}$. Since $V$ is transcendental,
so is $u$.  In (\ref{eq:U}), replace $U$ by $u$; then
$X=\ang{U}$, $v = \ang{U^{q-1}}$, and $z = \ang{U^{q^2-1}}$ are replaced by
$x=\ang u$, $V=\ang{u^{q-1}}$, and $D_{q+1}(V)=\ang{u^{q^2-1}}$. 
Eq.~(\ref{eq:v}) becomes
$VxF(x)=0$, so $F(x)=0$. This shows that every root
of $D_{q-1}(X)-V$ is a root of $F(X)$.
By Lemmas~\ref{lem:fwx} and~\ref{lem:distinctRoots},
$G$ divides $F$.

Let $H(X)=F(X)/G(X) \in \KK[X]$. Then $\deg_X(H)
=(q^2-1)-(q-1)^2=2(q-1)$.  Since $F(X)$ and $G(X)$ are monic,
so is $H(X)$. We claim $H=h$, where
$$h(X)=X^{2q-2}+VX^{q-1}+1.$$ 
Since $H$ and $h$ are monic of the same degree, it suffices to show
that the roots of $h$ in $\cj K$ are distinct and that $h(r)=0$ implies 
$G(r)\ne0$ and $F(r)=0$. 

If $h$ had a repeated root $r$, then $h(r)=h'(r)=0$. Here,
$h'(r)=Vr^{q-2}$ vanishes only at $r=0$, but $h(0)\ne 0$. 
Thus, $h$ has no repeated roots.

Next, $h(r)=0$ if and only if $V=(r^{2q-2}+1)/r^{q-1}=\ang{r^{q-1}}$. 
Since $V$ is transcendental, so is $r$. For each $w\in\F_q^\x$, 
$D_{q-1}(wr)-V=D_{q-1}(wr) - \ang{r^{q-1}}$, and the right side is nonzero
since $r^{q-1}\left(D_{q-1}(wx)-\ang{r^{q-1}}\right)$ is a nontrivial polynomial
in $\F_q[r]$ and $r$ is transcendental.  This shows $G(r)\ne0$. 

Finally, since $D_{q+1}(V)=D_{q+1}(\ang{r^{q-1}}) = \ang{r^{q^2-1}}$,
\begin{eqnarray*}
VF(r) &=& V\,r^{q^2-1}+D_{q+1}(V)\,r^{q-1} + V^q \\
&=& \ang{r^{q-1}}\,r^{q^2-1} + \ang{r^{q^2-1}}\,r^{q-1} + \ang{r^{q(q-1)}} \\
&=& (r^{q^2+q-2}+r^{q^2-q})+(r^{q^2+q-2}+r^{q-q^2})+(r^{q^2-q}+r^{q-q^2})\\
&=& 0.
\end{eqnarray*}
This shows $F(r)=0$. The roots of $h$ provide
$2q-2$ distinct roots of $H=F(X)/G(X)$, therefore
$h=H$ and $F(X)=h(X)G(X)$. 
This proves~(\ref{aczIdentity}). \qed

\end{document}